\documentclass[11pt,a4paper,twoside]{amsart}
\usepackage{amsmath}
\usepackage{amssymb}
\usepackage{latexsym}

\newcommand{\co}{\colon\thinspace}    
\newcommand{\fnote}[1]{\footnote{\small sharp1}}
\newcommand{\inv}{^{-1}}              

\newcommand{\N}{{\mathbb N}}
\newcommand{\Z}{{\mathbb Z}}
\newcommand{\R}{{\mathbb R}}

\newcommand{\T}{{\mathbb T}}

\newcommand{\spt}{\mbox{supp}}
\newcommand{\Azero}{\mathcal{A}_0}
\newcommand{\Ac}{\mathcal{A}_c}
\newcommand{\tildeAzero}{\tilde{\mathcal{A}_0}}

\newtheorem{theorem}{Theorem}
\newtheorem{proposition}[theorem]{Proposition}

\newtheorem{definition}[theorem]{Definition}
\newtheorem{lemma}[theorem]{Lemma}
\newtheorem{remark}[theorem]{Remark}

\title{Subsolutions of time-periodic Hamilton-Jacobi equations}
\author{Daniel Massart}
\date{\today}
\begin{document}

\begin{abstract}
We prove the existence of $C^{1}$ critical subsolutions of the Hamilton-Jacobi equation for a time-periodic Hamiltonian system. We draw a consequence for the Minimal Action functional of the system. 
\end{abstract}
\maketitle

\section{Introduction}
The purpose of this note is to generalize to time-periodic Hamiltonian systems some results that are known for autonomous (time-free) systems, namely Theorems 1.3 and 1.6 of \cite{FS}, and Theorem 1 of \cite{ijm}.

We call time-periodic Hamiltonian a $C^{2}$ function $H \co T^{\ast}M\times \T \longrightarrow \R$, where $M$ is a closed, connected manifold, and $\T$ is the unit circle, such that the restriction of $H$ to any subset $T^{\ast}_{x}M\times \left\{t\right\}$, for $(x,t) \in M\times \T$, is strictly convex and superlinear (see \cite{Mather91} which originated this line of research). We make the additional assumption that the Hamiltonian flow of $H$ is complete. The $\T$ factor is understood as a periodic dependance on time, whence the name. The Hamilton-Jacobi equation (HJc) is
	\[  \frac{\partial u}{\partial t}+H(x,\frac{\partial u}{\partial x},t) = c
\]
where the unknown $u$ is a $C^{1}$ function $M\times \T \longrightarrow \R$, and $c \in \R$ is a constant.
In general there may be no solution at all. One possible way around this fact is to look for solutions in a weak sense, say, viscosity solutions (see\cite{Fathi}, \cite{BeRo}). Another is to look for subsolutions, i.e.
 $C^{1}$ functions $u$ such that
 	\[  \frac{\partial u}{\partial t}+H(x,\frac{\partial u}{\partial x},t) \leq c. 
\]
The two approaches turn out to be connected, as shown by \cite{FS}. 

Since $M\times \T$ is compact, any function is a subsolution for a sufficiently large $c$, so the set $I$ of $c \in \R$ such that (HJc) has a subsolution is not empty. 
A subsolution of (HJc) is a subsolution of (HJc'), for any $c'\geq c$, so $I$ is an interval, unbounded to the right. Due to the convexity and superlinearity of $H$, and to the compactness of $M\times \T$, $H$ is bounded below, so  $I$ must be bounded to the left. Its infimum is called the critical value of $H$, and denoted $\alpha (H)$.

It is natural to ask whether $I$ is closed, i.e. whether $\alpha (H) \in I$. A subsolution of (HJ$\alpha (H)$), if it exists, is called critical. When $H$ is autonomous  the answer to the latter question is provided by  Theorem 1.2 of (\cite{FS} :
\begin{theorem}[Fathi-Siconolfi]
There exists a $C^{1}$ critical subsolution.
\end{theorem}
We extend this theorem to the time-periodic case in Section \ref{section_soussol}. The idea of  (the first step of) the  proof is borrowed from \cite{BB}, and uses the estimates of Section 2.2, which were proved in \cite{ijm} for the autonomous case.

The Hamiltonian $H$ being convex and superlinear, we may take advantage of the Lagrangian formulation of Classical Mechanics. Define
\[
	\begin{array}{rcl}
L\co   TM\times \T & \longrightarrow & \R \\
(x,v,t) & \longmapsto & \sup_{p \in T^{\ast}_{x}M}\left\{<p,v>-H(x,p,t)\right\}
\end{array}
\]
then $L$ is $C^{2}$, fiberwise strictly convex and superlinear. It defines, via the Euler-Lagrange equation, a flow $\Phi_t$ on $TM\times \T$ which is complete since it is the conjuguate, under Legendre Transform, of the Hamiltonian flow of $H$.

Define $\mathcal{M}_{inv}$
to be the set of $\Phi_t$-invariant, compactly supported, Borel probability measures on $TM\times \T$.
Mather showed that the function (called action of the Lagrangian on measures)
	\[
	\begin{array}{rcl}
\mathcal{M}_{inv} & \longrightarrow & \R \\
\mu & \longmapsto & \int_{TM\times \T}	Ld\mu
\end{array}
\]
is well defined and has a minimum.  It turns out that this minimum is $-\alpha(H)$. For this reason $\alpha(H)$ is also denoted $\alpha(L)$. A measure achieving the minimum is called $L$-minimizing.

One drawback of this characterization of the critical value is that when you want to test the minimality of a measure, you first need to check invariance. With this in mind,  an important corollary of Theorem 1.3 of \cite{FS} is  Theorem 1.6 of \cite{FS}, which is itself an elaboration on a theorem proved by Ma\~n\'e in \cite{Mane96}, and was proved by Bangert (\cite{Bangert99}) in the special case when the Lagrangian is a Riemannian metric.
\begin{definition}
A  probability measure $\mu$ on $TM$ is called closed if 
	\[\int_{TM} \left\|v\right\|d\mu (x,v)<+\infty,
\]
 and for every smooth function $f$ on $M$, we have 
\[
\int_{TM\times \T} df(x).v d\mu(x,v)=0.
\]
\end{definition}
Mather proved in \cite{Mather91} that every invariant measure is closed.
\begin{theorem}[Fathi-Siconolfi]\label{thm_fermees}
We have
	\[ -\alpha (0)= \min \left\{\int_{TM}	Ld\mu \co \mu \mbox{ is closed }\right\}.
\]
Moreover, every closed measure that achieves the minimum above is invariant under the Euler-Lagrange flow of $L$, and is thus a minimizing measure. 
\end{theorem}
The strength of this theorem is that it allows to work with measures without having to verify a priori that they are invariant. We give an appropriate definition of a closed measure for a time-periodic Hamiltonian system in Section \ref{closed_measures}, and indicate how the proof of Theorem \ref{thm_fermees} carries over to that case. 

The critical value is thus  a useful tool for selecting interesting invariant subsets ; for instance the supports of minimizing measures (Mather set), or the Aubry set (see below). The following classical trick gives us more milk from the same cow.  If $\omega$ is a closed one-form on  $M$, then $L-\omega$ is again a convex and superlinear Lagrangian, and it has the same Euler-lagrange flow as $L$. Besides, by Mather's Lemma (invariant measures are closed) if $\mu \in \mathcal{M}_{inv}$, the integral $\int_{TM\times \T}	\omega d\mu$ only depends on the cohomology class of $\omega$.  Then the minimum over $\mathcal{M}_{inv}$ of $\int (L-\omega)d\mu$ is actually a function of the cohomology class of $\omega$, the opposite of which is called the $\alpha$-function of the system. An  $(L-\omega)$-minimizing measure is also called $(L,\omega)$-minimizing or $(L,c)$-minimizing if $c$ is the cohomology of $\omega$. To sum up
	\[
	\begin{array}{rcl}
\alpha_{L} \co H^1 (M,\R) & \longrightarrow & \R \\
c & \longmapsto & 
-\min \left\{\int_{TM\times \T}	(L-\omega)d\mu \co 
\mu \in \mathcal{M}_{inv}\;
\left[\omega\right]=c\right\}.
\end{array}
\]
 In particular $\alpha (L)= \alpha_L (0)$. We shall omit the subscript $L$ when no ambiguity is possible. 
Mather proved that $\alpha$ is  convex and  superlinear. The analogy with the Lagrangian goes no further ; in general $\alpha$ is neither stricly convex, nor $C^{1}$ (see \cite{gafa}). The regions where $\alpha$ is not stricly convex (being convex, it must then be affine) are called faces of $\alpha$. By  Proposition 6 of \cite{ijm} (see \cite{Bernard_Fourier} for the time-periodic case) changing the cohomology class within a given face does not select any new dynamics. The presence of faces is often correlated with some rationality properties of homology classes (see \cite{ijm}, Corollary 3). Understanding this phenomenon is the motivation for Theorem 1 of \cite{ijm}, which we extend to the time-periodic case in the last section. The proof uses both the estimates of Section 2, and the existence of a $C^{1}$ subsolution, instead of Whitney's Extension Theorem as in \cite{ijm}. 
 
\textbf{Acknowledgements} :  I thank the referee for his careful reading and insightful advice. It is a pleasure to acknowledge the great hospitality of the CIMAT in Guanajuato, M\'exico while working on this project. 

\section{Preliminaries}
\subsection{Some Weak KAM theory}
In this section we briefly recall a few definition, referring the reader to the bibliography (\cite{Fathi}, and \cite{CIS} for the time-periodic case) for more information. 
Define, for all $n \in \N$, 
	\[ 
	\begin{array}{rcl} h_n \co \left(M\times \T \right)\times \left(M\times \T \right) & \longrightarrow & \R \\
\left((x,t),(y,s)\right)& \longmapsto & \min \int^{s+n}_{t}L(\gamma,\dot\gamma,t)dt	+n\alpha(0)
\end{array}
\]
where the minimum is taken over all absolutely continuous curves \\
$\gamma \co \left[t,s+n\right]\longrightarrow M$ such that $\gamma (t)=x$ and $\gamma (s+n)= y$. Note that we abuse notation, denoting by the same $t$ an element of $\T=\R/\Z$ or the corresponding point in $\left[0,1\right[$. The Peierls barrier is then defined as
	\[ \begin{array}{rcl} h \co \left(M\times \T \right)\times \left(M\times \T \right) & \longrightarrow & \R \\
\left((x,t),(y,s)\right)& \longmapsto & \liminf_{n \rightarrow \infty} h_n \left((x,t),(y,s)\right).
\end{array}
\]
The Aubry set is
	\[\Azero := \left\{(x,t) \in M\times \T  \co h\left((x,t),(x,t)\right) =0 \right\}.
\]
We say a function $f \co M\times \T \longrightarrow \R$ is $(L,\alpha(0))$-dominated if for every absolutely continuous curve $\gamma \co \left[a,b\right]\longrightarrow M$ with $b\geq a$ we have
	\[ \int^{b}_{a} \left(L(\gamma,\dot\gamma,t)+\alpha(0)\right)dt  \geq f\left(\gamma (b),b)\right)-f\left(\gamma (a),a)\right).
\]
Such  functions exist and are Lipschitz (\cite{Fathi}, Lemma 4.2.2), hence almost everywhere differentiable by Rademacher's theorem ; wherever the derivative exists, they are subsolutions of the Hamilton-Jacobi equation (see \cite{FS}), that is
	\[ \frac{\partial f}{\partial t}+H(x,\frac{\partial f}{\partial x},t) \leq \alpha(0).
\]

 A forward (resp. backward) weak KAM solution is a function $u$ which is $(L,\alpha(0))$-dominated and, for every $(x,t) \in M\times \T$, there exists an absolutely continuous curve 
 $\gamma \co \left[t, +\infty \right]\longrightarrow M$ (resp. 
 $\gamma \co \left[-\infty , t \right]\longrightarrow M$
 such that $\gamma (t)=x$ and, for every $s \in\left[t, +\infty \right]$ (resp. 
 $s \in\left[-\infty ,t \right]$), we have
	\[\int^{s}_{t} \left(L(\gamma(t),\dot\gamma(t),t)+\alpha(0)\right)dt = 
	u\left(\gamma(s), \dot\gamma(s),s \right)-u\left(\gamma(t), \dot\gamma(t),t\right).
\]
For every forward weak KAM solution $u_+$ there exists a unique backward weak KAM solution $u_-$ such that $u_+\leq u_-$,  $u_+=u_-$ in $\Azero$ (\cite{Fathi}, Theorem 5.12). The pair $(u_+,u_-)$ is then called a weak KAM conjuguate pair. It is a remarkable fact that for all $(x,t)$, $(y,s)$ in $M\times \T$
	\[ h \left((x,t), (y,s) \right)=\sup \left\{ u_- (y,s) -u_+(x,t)\right\}
\]
where the supremum is taken over all weak KAM conjuguate pairs $(u_+,u_-)$ (\cite{Fathi}, Corollary 5.37).
\subsection{An estimate}
To clear up the notation, we  assume $\alpha(0)=0$ by replacing $L$ with $L-\alpha(0)$ .
Take $\epsilon >0$. Let $N(\epsilon) \in \N^*$ be the smallest integer such that 
	\[
	\forall n \geq N(\epsilon) , \forall (x,\tau),(y,\sigma)\in M\times \T,
	h_n ((x,\tau),(y,\sigma)) \geq h ((x,\tau),(y,\sigma))-\epsilon.
\]
Let $(u_-,u_+)$  be a weak KAM conjugate pair such that $(u_- - u_+)\inv (0)= \Azero$.
Define $A_{\epsilon}:= (u_- - u_+)\inv ([2\epsilon,+\infty[)$.
Let $a,b$ be elements of $\R \cup \pm \infty$ and let $\gamma : [a,b]\longrightarrow M$ be an absolutely continuous curve. 
Denote by $\mbox{Leb}$ the normalized Lebesgue measure on $\R$, by $\mbox{Int}$ the integer part, and set $\mu_{\gamma}([a,b])= \mbox{Leb}(\gamma\inv (A_{\epsilon}))$. Then 
\begin{lemma}\label{formule}
We have :
	\[\int^{b}_{a}L(\gamma(t),\dot\gamma(t),t)dt \geq u_+(\gamma(b),b)-u_+(\gamma(a),a)+ 
\epsilon \mbox{Int}(\frac{\mu_{\gamma}([a,b])}{N(\epsilon)}).
\]
\end{lemma}

\proof

Define inductively a sequence in $ \R \cup \pm \infty$ by $t_0 := a$ and 
	\[ t_{i+1} := \max \left\{ t_i \leq t \leq b \co t-t_i \geq  N(\epsilon), \mu_{\gamma}([a,b])\leq N(\epsilon) \right\}
\]
Set $n_i := \mbox{Int}(t_{i+1}-t_i)$ ; we have $n_i \geq N(\epsilon)$.
 Note that $\forall i \geq 1, (\gamma (t_i),t_i) \in A_{\epsilon}$ ; this is the reason why we need a max in the above formula. Also, denoting $ n = \max \left\{ i \co t_i \leq b \right\}$, we have 
	\[ n = \mbox{Int}(\frac{\mu_{\gamma}([a,b])}{N(\epsilon)})
\]
since $\mu_{\gamma}([t_i,t_{i+1}])= N(\epsilon)$.

Now, we have
\begin{eqnarray*}
\int^{b}_{a}L(\gamma(t),\dot\gamma(t),t)dt  & = & \sum_{i=0}^{ n-1}  \int^{t_{i+1}}_{t_i}L(\gamma(t),\dot\gamma(t),t)dt 
+	\int^{b}_{t_n}L(\gamma(t),\dot\gamma(t),t)dt \\
& \geq & 
\sum_{i=0}^{ n-1}  h_{n_i}((\gamma (t_{i+1}),t_{i+1}),(\gamma (t_{i}),t_{i}))\\
&&
+u_+ (\gamma (b),b)-u_+ (\gamma (t_{n}),t_{n}).
\end{eqnarray*}
Since $n_i \geq  N(\epsilon)$, we have 
\begin{eqnarray*}
	h_{n_i}((\gamma (t_{i+1}),t_{i+1}),(\gamma (t_{i}),t_{i})) & \geq &
	h((\gamma (t_{i+1}),t_{i+1}),(\gamma (t_{i}),t_{i}))-\epsilon \\
& \geq & u_-(\gamma (t_{i+1}),t_{i+1})-	u_+(\gamma (t_{i}),t_{i})-\epsilon 
\end{eqnarray*}
whence 
\begin{eqnarray*}
\int^{b}_{a}L(\gamma(t),\dot\gamma(t),t)dt & \geq & 	
\sum_{i=0}^{ n-1} \left[ u_-(\gamma (t_{i+1}),t_{i+1})-	u_+(\gamma (t_{i}),t_{i})-\epsilon \right]\\
&&
+u_+ (\gamma (b),b)-u_+ (\gamma (t_{n}),t_{n})\\
&=& 
\sum_{i=1}^{ n-1} \left[ u_-(\gamma (t_{i+1}),t_{i+1})-	u_+(\gamma (t_{i}),t_{i})-\epsilon \right]\\
&&
+u_+ (\gamma (b),b)-u_+ (\gamma (a)),a)
\end{eqnarray*}
and, because $(\gamma (t_{i+1}),t_{i+1})\in A_{\epsilon}$, 
	\[ \int^{b}_{a}L(\gamma(t),\dot\gamma(t),t)dt \geq  u_+ (\gamma (b),b)-u_+ (\gamma (a),a) +n\epsilon	
\]
which proves the Lemma.
\qed
\subsection{Consequence of the estimate}

\begin{lemma}\label{chi}
There exists a $C^{2}$ non-negative function 
$W:  M\times \T \longrightarrow \R $
which is positive outside   $\Azero$ and zero inside  $\Azero$,
 such that $\alpha(L-W)=\alpha(L)$ and $\Azero (L-W)=\Azero(L)$.
\end{lemma}
\proof 
First we point out that, denoting $\chi_{\epsilon}$ the characteristic function of $A_{\epsilon}$, Lemma \ref{formule} may be rewritten 
	
\begin{equation}\label{chi_epsilon}
\int^{b}_{a}L(\gamma(t),\dot\gamma(t),t)dt \geq	 	u_+(\gamma(b),b)-u_+(\gamma(a),a)+ 
\frac{\epsilon}{ N(\epsilon)}\int^{b}_{a}\chi_{\epsilon}(\gamma(t),t)dt  
-\epsilon 
\end{equation}
since for each $i$ we have 
	\[\int^{t_{i+1}}_{t_{i}}\chi_{\epsilon}(\gamma(t),t)dt = N(\epsilon)
\]
and 
	\[\int^{b}_{t_{n}}\chi_{\epsilon}(\gamma(t),t)dt \leq N(\epsilon).
\]
The map 
	\[\chi := \sup_{n \in \N}\frac{2^{-n}}{ N(2^{-n})}\chi_{2^{-n}}
\]
is integrable by Lebesgue's Monotone Convergence Theorem. So, taking the supremum over $n \in \N$ in Equation (\ref{chi_epsilon}) we get
	
	\[\int^{b}_{a}L(\gamma(t),\dot\gamma(t),t)dt \geq
	u_+(\gamma(b),b)-u_+(\gamma(a),a)+ 
\int^{b}_{a}\chi(\gamma(t),t)dt  
- 1.
\]
Now pick a $C^{2}$  function 
$W:  M\times \T \longrightarrow \R $
which is positive outside of  $\Azero$, and such that 
	\[\forall (x,t) \in M\times \T, \; 0 \leq W(x,t) \leq \chi (x,t).
\]
	 
First let us verify that such a function exists. For every $n$ in $\N$ we can find a $C^{2}$ map $W_n:  M\times \T \longrightarrow \R $ with $C^{2}$-norm $\leq 1$ and such that
	\[
	\forall (x,t) \in M\times \T, \; 0 \leq W_n(x,t) \leq \frac{2^{-n-1}}{ N(2^{-n})}\chi_{2^{-n}} (x,t).
\]
	
 Now consider
$W:=\sum_{n\geq 0}W_n$, then $W$ is $C^{2}$, non-negative, and 
	\begin{eqnarray*}
	\forall (x,t) \in A_{n+1}\bigcap M\times \T \setminus A_{n},\;
	W(x,t)& \leq & \sum_{k\geq n+1}\frac{2^{-k-1}}{ N(2^{-k})} \\
	& \leq & \frac{2^{-n}}{ N(2^{-n})}\\
	& \leq & \chi (x,t)
\end{eqnarray*}
the latter inequality being true because $(x,t) \notin A_{n}$. 

It remains to be seen that $\alpha(L-W)=\alpha(L)$.

First, note that since $W$ is non-negative, for any real number $c$, a subsolution of (HJc) for $H+W$ is also a subsolution of (HJc) for $H$ so  $0=\alpha(H)\leq \alpha (H+W)$.

 Conversely, let $\mu$ be an ergodic $(L-W)$-minimizing measure and let $\gamma \colon\thinspace \R \longrightarrow M$ be a curve such that $(\gamma, \dot{\gamma}, t)$ is a $\mu$-generic orbit. We have, for all $s,t$ in $\R$ :
\begin{eqnarray*}
&&\int^{b}_{a}\left\{L(\gamma(t),\dot\gamma(t),t)-W(\gamma(t),t)\right\}dt 
 \geq \\
&& u_+(\gamma(b),b)-u_+(\gamma(a),a)+ 
\int^{b}_{a}\left\{\chi(\gamma(t),t)-W(\gamma(t),t)\right\}dt -1\\
&& \geq 
u_+(\gamma(b),b)-u_+(\gamma(a),a)-1	
\end{eqnarray*}
 thus by Birkhoff's Ergodic Theorem
  \[\int (L-W)d\mu \geq 0.
\]This proves that $0=\alpha(L)\geq \alpha (L-W)$ so $\alpha (L)= \alpha (L-W)$.

Let us pause for a moment to prove

\begin{proposition}\label{stricte}
There exists a critical subsolution which is strict at every point of $M\times \T ^{1}\setminus \Azero$.
\end{proposition}
\begin{remark}
The autonomous case of this Proposition (\cite{FS}, Proposition 6.1) is the first step of the proof of  Theorem 1.3 of \cite{FS}. 
The idea of the proof that follows is borrowed from \cite{BB}.
\end{remark}
\proof
Take a weak KAM solution $u$ for $L-W$, where $W$ is given by Lemma \ref{chi}. Recall that the Hamiltonian corresponding to $L-W$ under Legendre transform if $H+W$. At every point of differentiability of $u$ we have 
	\[ \frac{\partial u}{\partial t}+H(x,\frac{\partial u}{\partial x},t) +W(x,t) \leq \alpha_{L-W}(0)= \alpha_{L}(0)
\]
that is,
\[ \frac{\partial u}{\partial t}+H(x,\frac{\partial u}{\partial x},t) \leq  \alpha_{L}(0) -W(x,t)
\]
so $u$ is a  subsolution for $L$, strict outside of $\Azero$.
\qed

Observe that, since we know from \cite{CIS} that any critical subsolution is actually a  solution of (HJ$\alpha(H)$) in $\Azero$, the latter Proposition implies the following characterization of the Aubry set  :
\begin{proposition}\label{charac}
A point $(x,t) \in M\times \T$ is in $\Azero$ if and only if no critical subsolution of (HJ) is strict at $(x,t)$.
\end{proposition}
Now let us come back to the proof of Lemma \ref{chi}. We still have to find $W$ such that $\Azero(L-W)=\Azero(L)$. First note that since $W$ is non-negative, and $0=\alpha(H)\leq \alpha (H+W)$, any critical subsolution of (HJ) for $H+W$ is also a critical subsolution of (HJ) for $H$. Besides, $W$ being positive outside $\Azero$, such a subsolution is strict (for $H$) outside $\Azero$. By Proposition \ref{charac}, this implies $\Azero(L-W) \supset \Azero(L)$.

For the converse inclusion we may need to modify $W$. Assume there exists a $W_1$ such that $0\leq W \leq W_1$, all inequalities being strict outside $\Azero$, and 
$\alpha(H+W_1)= \alpha (H+W)$. This can be achieved by replacing $W$ with $W/2$ and taking $W$ as $W_1$.
Then a critical subsolution for $H+W_1$ is a also a critical subsolution for $H+W$, and it is strict for $H+W$ outside $\Azero$, which proves $\Azero(L-W) \subset \Azero(L)$.

\qed

\section{Subsolutions}\label{section_soussol}
Now we extend to the time-periodic case  Theorem 1.3 of \cite{FS} :
\begin{theorem}\label{soussolC1}
There exists a $C^{1}$ critical subsolution which is strict at every point of $M\times \T ^{1}\setminus \Azero$.
\end{theorem}

At this point we assume the reader has Theorem 9.2 of \cite{FS} before his eyes and explain how it applies. Take
\begin{itemize}
	\item $N:=M\times \T ^{1}$
	\item $f:=u$ given by  Proposition \ref{stricte}
	\item $A := \Azero (L) = \Azero (L-W)$
	\item $B:=$  the domain of $du$ ; $B$ has full measure and $du$ is defined in $B$ and continuous in $A$
	\item since we do not require the $C^{1}$ subsolution to approximate the strict subsolution, we do not need to specify $\epsilon$ 
  \item
	\[ F := \left\{(x,p,t,\tau) \in T\left( M\times \T ^{1}\right)\setminus \Azero \co
	\tau+H(x,p,t)\leq \alpha_{L} (0) - W(x,t)  \right\}
\]
\item
	\[ O := \left\{(x,p,t,\tau) \in T\left( M\times \T ^{1}\right)\setminus \Azero \co
	\tau+H(x,p,t)< \alpha_{L} (0) - \frac{1}{2}W(x,t)  \right\}.
\]
\end{itemize}
Then Theorem 9.2 of \cite{FS} yields a function $g$ that is the required $C^{1}$ critical subsolution, strict at every point of $M\times \T ^{1}\setminus \Azero$.
\qed
\subsection{Closed measures}\label{closed_measures}
 If we are going to extend Theorem \ref{thm_fermees} to time-periodic systems we have to integrate functions on $T(M\times \T)$ with respect to measures that are only defined on $TM\times \T$. The crucial point in the proof of Mather's lemma is that invariant measures are supported on curves in $TM$ of type $(\gamma(t),\dot{\gamma}(t))$. In the time-dependant setting we are considering curves in $M\times \T$ of type $(\gamma(t),t)$ so their velocities are $(\gamma(t),t,\dot{\gamma}(t),1)$. So the measures on $T(M\times \T)$ that we shall use are concentrated on the hypersurface 
 $\left\{(x,t,v,1) \co (x,v,t) \in TM\times \T \right\}$ in $T(M\times \T)$.
This leads to the following
\begin{definition}
A  probability measure $\mu$ on $TM\times \T$ is called closed if 
\[
\int_{TM\times \T} \left\|v\right\|d\mu (x,v,t) < + \infty,
\]
 and for every smooth function $f$ on $M\times \T$, we have 
\[
\int_{TM\times \T} df(x,t).(v,1) d\mu(x,v,t)=0.
\]
\end{definition}
Then Mather's lemma and its proof carry over without modification.

Let us sketch briefly how the proof of Theorem 1.6 of \cite{FS} applies to the time-periodic case. The first part of the proof consists of showing that a closed measure that realizes the minimum is supported inside $\Azero$. To make it work in the time-periodic case it suffices to replace every occurence of  $H(x,d_xu)$ by $\partial_t u +H(x,\partial_x u,t)$. Then  apply Proposition 10.3 of \cite{FS} with $N=M\times \T$ instead of $M$, and you're done.
\section{Minimal Action}
\subsection{Preliminaries}
Since the $\alpha$-function of $L$ is convex, at every point its graph has a supporting hyperplane. We call face of $\alpha$ the intersection of the graph of $\alpha$ with one of its supporting hyperplane. By Fenchel (a.k.a. convex) duality it is equivalent to study the differentiability of $\beta$  or to study the faces of $\alpha$. 
If $c$ is a cohomology class, we call $F_c$ the largest face of $\alpha$ containing $c$ in its relative interior, and $\mbox{Vect}F_c$ the underlying vector space of the affine space it generates in $H^{1}(M,\R)$. We call $\tilde{V}_c$ the underlying vector space of the affine space generated by  pairs $(c', \alpha (c')-\alpha (c))$ where $c' \in F_c$.
Replacing, if necessary, $L$ by $L-\omega$ where $[\omega]=c$, we only need consider the case when $c=0$.
Likewise, replacing $L$ with $L-\alpha(0)$ we may assume $\alpha(0)=0$.

\begin{definition} 
Let $\tilde{E}_0$ be the set of $(c,\tau) \in H^1 (M\times \T,\R)=H^1 (M,\R)\times H^1 ( \T,\R)$ such that there exists  a smooth closed one-form $\omega$ on $M\times \T$ with $[\omega] = (c,\tau)$ and $\spt (\omega)\cap \Azero = \emptyset$. Let $E_0$ be the canonical projection of $\tilde{E}_0$ to $H^1 (M,\R)$.
\end{definition}
\begin{definition}
Let $\tilde{G_0}$ be the set of $(c,\tau) \in H^1 (M\times \T,\R)=H^1 (M,\R)\times H^1 ( \T,\R)$ such that there exists a continuous closed one-form $\omega$ on $M\times \T$ with $[\omega] = (c,\tau)$ and 
	\[\omega (x,t,v,\tau)=0 \;\forall (x,t) \in \Azero \subset M\times \T, \; \forall (v,\tau)\in 
	T_{(x,t)}M\times \T.
\]
 Let $G_0$ be the canonical projection of $\tilde{G}_0$ to $H^1 (M,\R)$.
\end{definition}
Now we can state the main result of this section 
\begin{theorem}\label{inclusion} The following inclusions hold true : 
	\[E_0 \subset \mbox{Vect}F_0 \subset G_0.
\]
\end{theorem}
In view of the above definitions we shall need to integrate one forms on $M\times \T$ with respect to invariant measures. We denote by
$\int \omega d\mu$
the expression
$$
\int_{TM\times \T}\omega_{(x,t)}\cdot (v,1)d\mu(x,v,t).
$$
The following lemma is useful.
\begin{lemma}\label{1forme}
If $\omega$ is a closed one form on $M\times \T$, with 
$[\omega]=(c,\tau) \in H^{1}(M,\R)\times H^{1}(\T,\R)$, and $\mu$ is an $(L,c)$-minimizing measure, then
\[
\int (L-\omega)d\mu = -\alpha (c) -\tau.
\]
\end{lemma}
\proof
Consider a closed one-form $\omega_1$  on $M$ such that $\left[\omega_1\right]=c$. Denote $\tilde{\tau}$ the constant one-form $\tau dt$ on $\T$. Then $\omega_1\oplus \tilde{\tau}$ is a one-form on $M\times \T$, cohomologous to 
$\omega$. Let $f$  be a smooth function on $M\times \T$ such that $(\omega_1,\tilde{\tau})=\omega +df$. Then by Mather's lemma (invariant measures are closed) 
	\[\int (L-\omega)d\mu = \int (L-(\omega_1\oplus \tilde{\tau}))d\mu. 
\]
On one hand $\int (L-\omega_1)d\mu = -\alpha(c)$ since $\mu$ is $(L,c)$-minimizing. On the other hand, since $\mu$ is a probability measure, we have 
$\int \tilde{\tau} d\mu =  \int \tau d\mu =\tau$. The lemma is proved.
\qed

\subsection{Proof of $E_0 \subset \mbox{Vect}F_0$}
Pick $ c \in E_0$. Let $\tau \in H^{1}(\T,\R)$ and $\omega$ a closed one-form on $M\times \T ^{1}$ 
 be such that $\spt (\omega) \cap \Azero = \emptyset $ and $\left[\omega \right]=(c,\tau)$.
 Since  $\spt (\omega)$ is compact there exists $\epsilon>0$ such that 
	\[u_- (x,t) - u_+(x,t) \geq 2\epsilon \; \forall (x,t) \in \spt (\omega).
\]
By {\em a priori} compacity there exists a compact subset $K$ in $TM\times \T ^{1}$ such that for all $\theta\in \left[-1,1 \right]$, for all $L+\theta\omega$-minimizing measure $\mu$, the support of $\mu$ is contained in $K$. Let $\delta$ be such that
	\[ \forall (x,v,t) \in K, |\delta\omega_{(x,t)}(v,1)| \leq \frac{\epsilon}{N(\epsilon)}.
\]
   Let $\mu$ be an ergodic  $(L+\delta\omega)$-minimizing measure and let $\gamma \co \R \longrightarrow M$ be a $\mu$-generic orbit. We have, for all $s\leq t$ : 
	\[ \left| \int^{s}_{t}\delta\omega (\gamma,\dot\gamma,t)dt \right| \leq (t-s)\frac{\epsilon}{N(\epsilon)} \leq
	\epsilon \mbox{Int}\left(\frac{t-s}{N(\epsilon)}\right) +\epsilon
\]
whence
	\[\int^{s}_{t}(L+\delta\omega) (\gamma,\dot\gamma,t)dt \geq u_+(\gamma(t),t)-u_+(\gamma(s),s)-\epsilon
\]

  thus by Birkhoff's Ergodic Theorem $\int (L+\delta\omega)d\mu \geq 0$.
Now by Lemma \ref{1forme} 
\[
\int (L+\delta\omega)d\mu = -\alpha (\delta c)-\delta\tau \mbox{ so }
\alpha (\delta c)\leq -\delta\tau.
\]
 Likewise, $\alpha (-\delta c)\leq \delta\tau$ thus 
$\alpha (\delta c) +\alpha (-\delta c)\leq 0$. On the other hand by convexity of $\alpha$ the reverse inequality is true : $\alpha (\delta c) +\alpha (-\delta c)\geq 0$ so the inequalities $\alpha (\delta c)\leq -\delta\tau$ and $\alpha (-\delta c)\leq \delta\tau$ are actually equalities. This means that $\alpha$ restricted to the line segment $\left[-\delta c,\delta c\right]$, is affine with slope $-\tau$, which proves that $-\tau=\alpha(c)$ and $\delta c\in  F_0$ whence $c\in  \mbox{Vect}F_0$. \qed

 \subsection{Proof of $\mbox{Vect}F_0 \subset G_0$}
 Pick $c$ in the interior of $F_0$. Note that by Proposition 6 of \cite{ijm} we have $\Ac=\Azero$.
 Take $\omega$ a smooth closed one-form on $M$ such that $[\omega]=c$.
 Let $u_0$ (resp. $u_1$) be a $C^{1}$ subsolution for $L$ (resp. $L-\omega$). Then for all $(x,v,t) \in \tildeAzero$, we have

\begin{eqnarray*}
	\frac{\partial u_0}{\partial x}(x,t)  & = & \frac{\partial L}{\partial v}(x,v,t)\\
		\frac{\partial u_1}{\partial x}(x,t)  & = & \frac{\partial L}{\partial v}(x,v,t)-\omega_x (v)
\end{eqnarray*}
 Observe that the Hamiltonian paired by Legendre transform with $L-\omega$ is $(x,p,t)\longmapsto H(x,p+\omega_x,t) := H_{\omega}(x,p,t)$. Thus
	\[\forall (x,t) \in \Azero \; \; 
	H_{\omega}(x,\frac{\partial u_1}{\partial x}(x,t),t)=H(x,\frac{\partial u_0}{\partial x}(x,t),t).
\]
On the other hand in $\Azero$ $u_0$ and $u_1$ are solutions of the Hamilton-Jacobi equation : 
\begin{eqnarray*}
\frac{\partial u_0}{\partial t}(x,t)	+ H(x,\frac{\partial u_0}{\partial x}(x,t),t) & = & \alpha (0) \\
\frac{\partial u_1}{\partial t}(x,t)	+ H_{\omega}(x,\frac{\partial u_1}{\partial x}(x,t),t) & = & \alpha (c)
\end{eqnarray*}
whence 
	\[\frac{\partial (u_1-u_0)}{\partial t}(x,t)=\alpha (c)-\alpha(0)\; \forall (x,v,t) \in \tildeAzero.
\]
 Consider the closed one-form $\tilde{\omega}$ on $M\times \T$ defined by 
	\[\tilde{\omega}_{(x,t)}(v,\tau):= \omega_{x}(v)+(\alpha (0)-\alpha(c))\tau.
\]
 The cohomology class of $\tilde{\omega}$ is $(c,\alpha (0)-\alpha(c))$ and $\tilde{\omega}=d(u_0-u_1)$ in $\Azero$ so replacing $\tilde{\omega}$ by the continuous one-form $\tilde{\omega}-d(u_0-u_1)$ we see that $c\in G_0$. \qed

{\small

\bigskip

\noindent

Math\'ematiques, Universit\'e Montpellier II, France\\
e-mail : massart@math.univ-montp2.fr
}

\end{document}